\documentclass[11pt,a4paper]{article}
\usepackage[numbers,sort&compress]{natbib}
\usepackage[colorlinks,
            linkcolor=blue,
            anchorcolor=green,
            citecolor=magenta
            ]{hyperref}
\usepackage{mathrsfs,amssymb,amsfonts}
\usepackage{graphicx}
\usepackage{amsmath,amssymb,latexsym,color}
\usepackage[mathscr]{eucal}
\usepackage[numbers]{natbib}
\usepackage{bm}
\usepackage[top=1in, bottom=1in, left=1.25in, right=1.25in]{geometry}
\usepackage{tikz}
\usepackage{verbatim}
\usetikzlibrary{arrows,decorations.pathmorphing,backgrounds,positioning,fit,petri}
\usepackage{setspace}
\makeatletter

\newcommand{\Rmnum}[1]{\expandafter\@slowromancap\romannumeral #1@}
\makeatother
\newtheorem{theorem}{Theorem}[section]

\newtheorem{lemma}[theorem]{Lemma}
\newtheorem{remark}[theorem]{Remark}

\newtheorem{corollary}[theorem]{Corollary}
\newtheorem{example}[theorem]{Example}

\newcommand{\qed}{\hfill\Box\medskip}
\usepackage{geometry}
\usepackage{diagbox}
\begin{document}
\begin{spacing}{1.5}


\title{\bf  A characterization of graphs attaining a lower bound on the dissociation number \footnote{This research is supported by the National Natural Science Foundation of China (12301453), the Project of Scientific Research Fund of Hunan Provincial Science and Technology Department (2018WK4006).}
}

\author{Pingshan Li\footnote{ E-mail address: lips@xtu.edu.cn}
\\
{\footnotesize \em School of Mathematics and Computational Science, Xiangtan University,
Xiangtan, Hunan 411105, PR China}\\
{\footnotesize Key Laboratory of Intelligent Computing$\And$Information Processing of Education}\\
{\footnotesize Key Laboratory for Computation and Simulation in Science and Engineering}\\
{\footnotesize National Center for Applied Mathematics in Hunan}\\
}

 \date{}
 \date{}
 \maketitle

\begin{abstract}
The dissociation number ${\rm diss}(G)$ of a graph $G$  is defined as the maximum cardinality of a vertex subset $D\subseteq V(G)$ such that the induced subgraph $G[D]$ has maximum degree at most $1$. In 2023,  Bock, Pardey, Penso and Rautenbach [J. Graph Theory 103 (2023), 661--673] showed that ${\rm diss}(G)\ge n-\frac{1}{3}(m+k+c_1)$ where $n, m, k, c_1$ denote the numbers of vertices, edges, components and induced cycles of length $1$ modulo $3$ in $G$ and they also established a complete characterization of extremal graphs under the vertex-disjoint cycle constraint. Although this partial characterization offers very valuable structural insights, it excludes many simple extremal graphs. So they  posed a natural question: How can all extremal graphs be comprehensively characterized without restrictions? In this paper, we give a recursive characterization of all graphs attaining equality in this bound.
\medskip

\noindent {\em Key words:}  Dissociation sets; dissociation number; maximum; extremal graphs.

\medskip
\end{abstract}

\section{Introduction}

In this paper, we focus on simple, finite, and undirected graphs. A vertex-set $D$ of a graph $G$ is a \textit{dissociation set} if the induced subgraph $G[D]$ has maximum degree at most $1$. A dissociation set is called \textit{maximum} if it has maximum cardinality. The \textit{dissociation number} ${\rm diss}(G)$ of $G$ is the cardinality of a maximum dissociation set  in $G$. The problem of computing ${\rm diss}(G)$ was introduced by Yannakakis and was shown to be NP-hard for the class of bipartite graphs \cite{yannakakis}. Moreover, the problem remains NP-hard for planar line graphs of planar bipartite graphs with maximum degree $4$ \cite{orlovich}.
This complexity motivates the study of bounds for the dissociation number, and several tight lower bounds have been proposed to estimate its minimum possible value in different graph classes: if $G$ is a graph with $n$ vertices and $m$ edges, then

\begin{equation}
\begin{split}
    {\rm diss}(G)\ge
    \begin{cases}
    \frac{2n}{3} - \frac{m}{6}, &\text{\cite{BRESAR}},\\
      \frac{2n}{3}, &\text {if } G \text { is a tree \cite{BRESAR}},\\
      \frac{4}{3} \sum\limits_{u \in V(G)} \frac{1}{d_G(u) + 1}, & \text{if } G \text{ has no isolated vertex  \cite{BRESAR}},\\
    \frac{n}{\left\lceil \frac{\Delta + 1}{2} \right\rceil}, & \text{if } G \text{ has maximum degree } \Delta~\text{\cite{BRESAR}},\\
    \frac{n}{2}, &\text {if } G \text{ is outerplanar \cite{BRESAR}},\\
    \frac{2n}{r+2} - \frac{m}{(r+1)(r+2)}, &\text {if } m>0,\ r = \left\lceil \frac{m}{n} \right\rceil - 1 ~\text{\cite{BRESAR2}},\\
     \sum\limits_{u \in V(G)} \frac{1}{d_G(u) + 1} + \sum\limits_{uv \in E(G)} \binom{|N_G[u] \cup N_G[v]|}{2}^{-1}, & \text{\cite{6}}
    \end{cases}
\end{split}
\label{eq:result}
\end{equation}

Notably, Bock, Pardey, Penso and Rautenbach generalized the lower bound $\frac{2}{3}n$ for trees of order $n$ (as stated in (1)), and simultaneously strengthened the general lower bound $\frac{2}{3}n-\frac{m}{6}$ for a broad class of graphs. Their result is formalized as follows:

\begin{theorem}{\rm \cite{bock}}\label{theorem1} If $G$ is a graph with $n$ vertices, $m$ edges, $k$ components, and $c_1$ induced cycles of length $1$ modulo $3$, then
  \begin{equation}\label{(1)}
  {\rm diss}(G)\ge n-\frac{1}{3}(m+k+c_1). \tag{2}
  \end{equation}
\end{theorem}

As a crucial contribution of their paper, they established a characterization of the extremal graphs for Theorem \ref{theorem1} under the constraint that every two cycles are vertex-disjoint (see Theorem \ref{theorem2}(1)). To elaborate on the structure of these extremal graphs, we begin by introducing several foundational concepts.

 \begin{figure}[!htbp]
  \centering
  \includegraphics[width=0.85\textwidth]{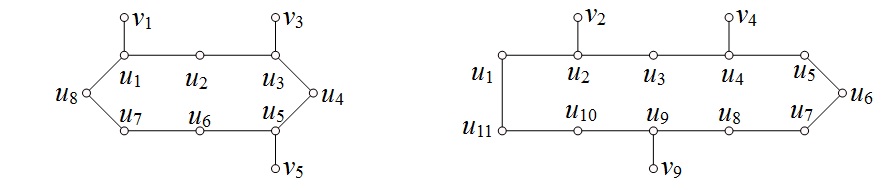}\\
  \caption{Two very good spiked cycles $C_8^{1,3,5}, C_{11}^{2, 4, 9}$}\label{very}
\end{figure}
The set $\{1, 2, \cdots, n\}$ is denoted as $[n]$. Let $\ell\ge 3$ and $1\le i_1<i_2<\cdots<i_k\le \ell$. A spiked cycle $C^{i_1, i_2, \cdots, i_k}_{\ell}$ is the graph obtained from the cycle $C=u_1u_2\cdots u_{\ell}u_1$ by adding $k$ new vertices $v_{i_j}$ and $k$ edges $u_{i_j}v_{i_j}$ for $j\in [k]$. The spiked cycle $C^{i_1, i_2, \cdots, i_k}_{\ell}$ is called \emph{very good} (see Figure \ref{very} for illustration) if
\begin{center}
\begin{itemize}
  \item $k\ge 2, \ell \not\equiv 1$ mod $3$,
  \item $i_{j+1}-i_j\equiv 2$ mod $3$ for every $j\in [k-1]$, and
  \item $\ell+i_1-i_k\equiv 1$ mod $3$.
\end{itemize}
\end{center}

\begin{figure}[!htbp]
  \centering
  \includegraphics[width=1\textwidth]{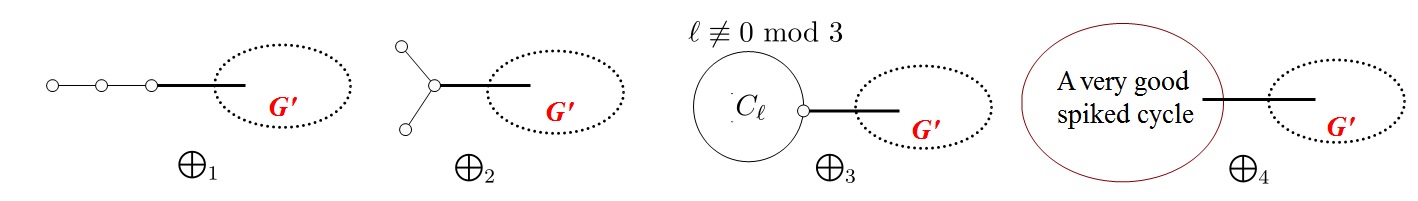}\\
  \caption{Operations constructing an extremal graph from a smaller extremal graph $G'$}\label{figure1}
\end{figure}

The operations $\bigoplus_1, \bigoplus_2, \bigoplus_3$ and $\bigoplus_4$ (see Figure \ref{figure1}) are bridge-adding operations between a given graph $G'$ and one of the following: a path $P_3$, a cycle $C_{\ell}$ with $\ell\not\equiv 0$ mod $3$ or  a very good spiked cycle.
The resulting graphs are denoted as $G'\bigoplus_1 P_3, G'\bigoplus_2 P_3, G'\bigoplus_3 C_{\ell}$ and $G'\bigoplus_4 C^*$ respectively. Let $\mathcal{C}$ be the set of all graphs that arise from the graphs in
\begin{center}
$\mathcal{C}_0=\{P_3\}\cup\{C_{\ell}: \ell\ge3~{\rm ~is~an~integer~and}~\ell \not\equiv 0~{\rm mod~} 3\}\cup\{C^*: C^*~{\rm is~a~very~good~spiked~cycle}\}$
\end{center}
 by repeated applications of $\bigoplus_i$ for $i\in [4]$.

\begin{theorem}{\rm\cite{bock}}\label{theorem2} {\rm (1)} A connected cycle-disjoint graph satisfies {\rm (\ref{(1)})} with equality if and only if it belongs to $\mathcal{C}$.

{\rm (2)} For every vertex $u$ of a graph $G$ in $\mathcal{C}$, $G$ has a maximum dissociation set not containing $u$.
\end{theorem}

\begin{figure}[!htbp]
  \centering
  \includegraphics[width=0.4\textwidth]{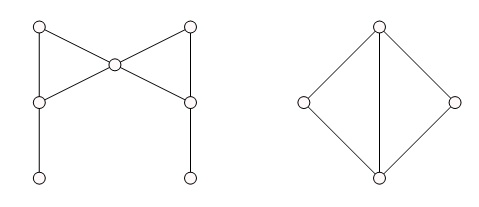}\\
  \caption{Two extremal graphs in Theorem \ref{theorem1} violate the vertex-disjoint cycle condition.}\label{figure+}
\end{figure}

Note that the restriction that every two induced cycles must be vertex-disjoint is so restrictive that many simple extremal graphs of Theorem \ref{theorem1} (as illustrated in Figure \ref{figure+}) are excluded. Bock et al. posed a natural question: How can all extremal graphs of Theorem \ref{theorem1} be comprehensively characterized without restrictions? This paper establishes a complete characterization of these extremal graphs (Theorem \ref{theorem3}). As preparatory groundwork, we also introduce two operations.


The operation $\bigoplus_5$ adds a new vertex $x$ and the three edges $xu,xv,xw$ to a connected graph $G_1$, where $u,v,w$ are distinct vertices of $G_1$. Put $S=\{u,v,w\}$. The operation is permitted if
\begin{enumerate}
  \item for every maximum dissociation set $D$ of $G_1$, either $|D\cap S|\ge2$, or $D\cap S=\{z\}$ and $d_{G_1[D]}(z)=1$;
  \item there is no induced path $P$ in $G_1$ of length $2$ modulo $3$ whose endpoints $s,t$ belong to $S$ and for which $V(P)\cap S=\{s,t\}$.
\end{enumerate}
Here and throughout the paper, the length of a path is its number of edges.

\begin{figure}[!htbp]
  \centering
  \includegraphics[width=0.6\textwidth]{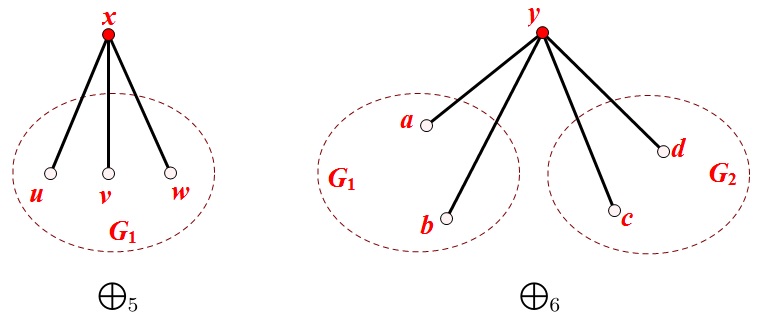}\\
  \caption{Operations constructing an extremal graph from smaller extremal graphs}\label{figure56}
\end{figure}

The operation $\bigoplus_6$ adds a new vertex $y$ and the four edges $ya,yb,yc,yd$ to the disjoint union of two connected graphs $G_1,G_2$, where $a,b$ are distinct vertices of $G_1$ and $c,d$ are distinct vertices of $G_2$. Put $S=\{a,b,c,d\}$. The operation is permitted if
\begin{enumerate}
  \item for every maximum dissociation set $D$ of $G_1\cup G_2$, either $|D\cap S|\ge2$, or $D\cap S=\{z\}$ and $d_{(G_1\cup G_2)[D]}(z)=1$;
  \item there is no induced $a$--$b$ path in $G_1$ or induced $c$--$d$ path in $G_2$ of length $2$ modulo $3$.
\end{enumerate}
The resulting graphs are denoted by $G_1\bigoplus_5\{x\}$ and $G_1\bigoplus_6G_2$, respectively; see Figure~\ref{figure56}.

Let $G_1,G_2$ be two disjoint connected graphs. We use $G_1\biguplus G_2$ to denote a graph obtained by adding an edge between an arbitrary vertex of $G_1$ and an arbitrary vertex of $G_2$. Let $\mathcal{G}$ be the smallest class of connected graphs containing $\mathcal{C}_0$ and closed under the operations $\bigoplus_i$ for $i\in[6]$ and under this bridge-joining operation.

\begin{theorem}\label{theorem3} A connected graph $G$ satisfies {\rm (\ref{(1)})} with equality if and only if $G\in\mathcal{G}$.
\end{theorem}

For convenience, we always use $n_i, m_i$ and $c_{1i}$ to denote $|V(G_i)|, |E(G_i)|$ and the number of induced cycles with length $1$ modulo $3$ in $G_i$, respectively. By Theorem \ref{theorem3}, if $G_1, G_2\in \mathcal{G}$, then

$\begin{array}{rl}
{\rm diss}(G_1\bigcup G_2)&= {\rm diss}(G_1)+{\rm diss}(G_2)\\
  & = (n_1-\frac{1}{3}(m_1+1+c_{11}))+(n_2-\frac{1}{3}(m_2+1+c_{12}))\\
& = (n_1+n_2)-\frac{1}{3}((m_1+m_2+1)+1+(c_{11}+c_{12})) \\
&={\rm diss}(G_1\biguplus G_2).
\end{array}$
\\Thus, we have the following corollary.
\begin{corollary}
A graph $G$ satisfies {\rm (\ref{(1)})} with equality if and only if every component of $G$ belongs to $\mathcal{G}$.
\end{corollary}

\section{The proof of the sufficiency of Theorem \ref{theorem3}}

We prove the sufficiency of Theorem \ref{theorem3} by proving a bit stronger result.

\begin{lemma}\label{lemma2} All graphs in $\mathcal{G}$ satisfy $(\ref{(1)})$ with equality. Furthermore, for every vertex $u$ of a graph $G$ in $\mathcal{G}$, $G$ has a maximum dissociation set not containing $u$.
\end{lemma}

{\noindent \bf Proof:} We prove both assertions by induction on the number $j$ of operations used in a construction from $\mathcal{C}_0$. By Theorem~\ref{theorem2}(1)--(2), the result holds for $j=0$. Suppose the assertions hold for all constructions using fewer than $j$ operations. We distinguish the following cases according to the final operation.

{\noindent\bf Case 1:} The final operation joins two graphs by a bridge.

Write $G=G_1\biguplus G_2$, where $G_1, G_2\in\mathcal{G}$ both have constructions using fewer than $j$ operations. This case includes the operations $\oplus_i$ for $i\in [4]$. Restricting a dissociation set of $G$ to the two input graphs gives
\begin{center}
 ${\rm diss}(G)={\rm diss}(G_1\biguplus G_2)\le {\rm diss}(G_1)+{\rm diss}(G_2)$.
\end{center}
By induction, ${\rm diss}(G_i)=n_i-\frac{1}{3}(m_i+1+c_{1i})$ for $i\in [2]$. By Theorem \ref{theorem1},
\begin{center}
$\begin{array}{rl}
{\rm diss}(G_1)+{\rm diss}(G_2)& = (n_1-\frac{1}{3}(m_1+1+c_{11}))+(n_2-\frac{1}{3}(m_2+1+c_{12}))\\
& = (n_1+n_2)-\frac{1}{3}((m_1+m_2+1)+1+(c_{11}+c_{12})) \\
&=n-\frac{1}{3}(m+1+c_1)\\
&\le{\rm diss}(G_1\biguplus G_2)\\
&={\rm diss}(G).
\end{array}$
\end{center}
By combining the two inequalities above, we can obtain ${\rm diss}(G)=n-\frac{1}{3}(m+1+c_1)$.

Suppose $u_1u_2$ is the edge between $G_1$ and $G_2$ where $u_i\in V(G_i)$. For an arbitrary vertex $x_1\in V(G_1)$, there exists a maximum dissociation set $D_1$ of $G_1$ that avoids $x_1$ by induction. Let $D_2$ be a maximum dissociation set of $G_2$ that does not contain $u_2$.
Then, $D_1\cup D_2$ is a maximum dissociation set of $G$ that does not contain $x_1$. Similarly, there exists a maximum dissociation set of $G$ avoiding any specified vertex of $G_2$.

{\noindent\bf Case 2:} The final operation is $\bigoplus_5$.

In this case, let $G=G_1\bigoplus_5\{x\}$, and let $D_1$ be an arbitrary maximum dissociation set of $G_1$.
Clearly, $D_1$ is a dissociation set of $G$. It follows that ${\rm diss}(G)\ge {\rm diss}(G_1)$.
Suppose that ${\rm diss}(G)>{\rm diss}(G_1)$. Since only one vertex is added, ${\rm diss}(G)={\rm diss}(G_1)+1$. Every maximum dissociation set $D$ of $G$ then contains $x$, and $D_0=D\setminus\{x\}$ is a maximum dissociation set of $G_1$. Put $S=N_G(x)$. Since $D$ is a dissociation set of $G$, either $D_0\cap S=\emptyset$, or $D_0\cap S=\{z\}$ and $z$ is isolated in $G_1[D_0]$. Both alternatives contradict condition (1) of $\bigoplus_5$. Hence ${\rm diss}(G)={\rm diss}(G_1)$.
Thus, every maximum dissociation set of $G_1$ is also a maximum dissociation set of $G$.
This means that there exists a maximum dissociation set of $G$ avoiding $x$. By induction, there exists a maximum dissociation set of $G_1$ avoiding any specified vertex of $G_1$. It follows that there exists a maximum dissociation set of $G$ avoiding any specified vertex of $G_1$.
To complete the proof of this case, we just need to show that ${\rm diss}(G)=n-\frac{1}{3}(m+1+c_1)$.

Every new induced cycle contains $x$. Deleting $x$ from such a cycle leaves an induced path whose endpoints belong to $N_G(x)$ and whose internal vertices avoid $N_G(x)$. By condition (2) of $\bigoplus_5$, this path cannot have length $2$ modulo $3$. Thus no new induced cycle has length $1$ modulo $3$, and the old induced cycles are unchanged. Consequently, $c_{11}=c_1$. By induction,
\begin{center}

$\begin{array}{rl}
  {\rm diss}(G_1)& =n_1-\frac{1}{3}(m_1+1+c_{11}) \\
   & =(n-1)-\frac{1}{3}((m-3)+1+c_{11})\\
    & =n-\frac{1}{3}(m+1+c_1).
\end{array}
$
\end{center}
Recall that ${\rm diss}(G_1)={\rm diss}(G)$. The result holds.

{\noindent\bf Case 3:} The final operation is $\bigoplus_6$.

Let $G_1\cup G_2$ be the graph before the last operation, and let $D_i$ be an arbitrary maximum dissociation set of $G_i$ for $i\in [2]$.
Clearly, $D_1\cup D_2$ is a dissociation set of $G$. It follows that ${\rm diss}(G)\ge {\rm diss}(G_1)+{\rm diss}(G_2)$. Arguing as in Case 2, we obtain that
 \begin{center}
 ${\rm diss}(G)={\rm diss}(G_1)+{\rm diss}(G_2)$.
\end{center}
Thus, every maximum dissociation set of $G_1\cup G_2$ is also a maximum dissociation set of $G$. This means there is a maximum dissociation set of $G$ avoiding the added vertex of operation $\bigoplus_6$. Note that there exists a maximum dissociation set of $G_i$ avoiding any specified vertex of $G_i$ by induction. Taking the union of suitably chosen maximum dissociation sets of \(G_1\) and \(G_2\) yields a maximum dissociation set of \(G\) avoiding any prescribed vertex of \(G_1 \cup G_2\).
 To complete the proof of this case, we just need to show that ${\rm diss}(G)=n-\frac{1}{3}(m+1+c_1)$.

Operation $\bigoplus_6$ reduces the number of components by $1$. Every new induced cycle consists of the added vertex together with an induced $a$--$b$ path in $G_1$ or an induced $c$--$d$ path in $G_2$. Condition (2) excludes new induced cycles of length $1$ modulo $3$. Thus, $c_{11}+c_{12}=c_1$. By induction,
\begin{center}
 $\begin{array}{rl}
                  {\rm diss}(G_1)+{\rm diss}(G_2) & =(n_1-\frac{1}{3}(m_1+1+c_{11}))+(n_2-\frac{1}{3}(m_2+1+c_{12})) \\
                   & =(n_1+n_2)-\frac{1}{3}(m_1+1+c_{11}+m_2+1+c_{12})\\
                   & =(n_1+n_2+1)-\frac{1}{3}((m_1+m_2+4)+1+(c_{11}+c_{12}))\\
                   &=n-\frac{1}{3}(m+1+c_1).
                \end{array}$
\end{center}
Recall that ${\rm diss}(G)={\rm diss}(G_1)+{\rm diss}(G_2)$. The result holds.
$\qed$

\section{The proof of the necessity of Theorem \ref{theorem3}}

\hspace{2.0em}Let $G$ be a connected graph satisfying (\ref{(1)}) with equality. By induction on its order $n$, we show that $G\in\mathcal{G}$. If $G$ is cycle-disjoint, then the result holds by Theorem~\ref{theorem2}(1). Thus, we may assume that $G$ is not cycle-disjoint.

Choose a vertex $u$ such that $G-u$ has at most $d_G(u)-2$ components. Such a vertex exists: if a block $B$ with at least three vertices is not a cycle, choose $u\in V(B)$ with $d_B(u)\ge3$. Since $B-u$ is connected, at least three neighbors of $u$ belong to the same component of $G-u$. Otherwise, all blocks with at least three vertices are cycles, and two cycle blocks share a vertex $u$; each of their two disjoint pairs of neighbors remains connected after deleting $u$. In either case, every component of $G-u$ contains a neighbor of $u$, giving the required bound.

Let $C_1,\ldots,C_t$ be the components of $G-u$, and let $c'_1$ be the number of induced cycles of length $1$ modulo $3$ in $G-u$. Thus $t\le d_G(u)-2$ and $c'_1\le c_1$.

{\noindent\bf Claim 1:} $t=d_G(u)-2$, $c'_1=c_1$, and ${\rm diss}(G-u)={\rm diss}(G)$. Moreover, every $C_i$ satisfies (\ref{(1)}) with equality.

{\noindent\bf Proof of Claim 1:} By Theorem~\ref{theorem1},
\begin{align*}
 n-\frac{m+1+c_1}{3}
 = {\rm diss}(G)
 &\ge {\rm diss}(G-u)\\
 &\ge (n-1)-\frac{(m-d_G(u))+t+c'_1}{3}\\
 &=n-\frac{m+1+c_1}{3}
   +\frac{d_G(u)-t-2+c_1-c'_1}{3}\\
 &\ge n-\frac{m+1+c_1}{3}.
\end{align*}
All inequalities are equalities. Since both $d_G(u)-t-2$ and $c_1-c'_1$ are nonnegative, they both vanish. Hence $t=d_G(u)-2$, $c'_1=c_1$, and ${\rm diss}(G-u)={\rm diss}(G)$. The bound in Theorem~\ref{theorem1} is additive over components, and each component satisfies it. Equality for $G-u$ therefore implies equality for every $C_i$. $\qed$

By induction, $C_i\in\mathcal{G}$ for every $i\in[t]$. In particular, Lemma~\ref{lemma2} provides a maximum dissociation set of each $C_i$ avoiding any prescribed vertex. Put $r_i=|N_G(u)\cap V(C_i)|$. Since $G$ is connected, $r_i\ge1$, and Claim 1 gives
\[
 \sum_{i=1}^{t}(r_i-1)=d_G(u)-t=2.
\]
Thus, after reordering the components, either two of the $r_i$ are $2$ and all others are $1$, or one of the $r_i$ is $3$ and all others are $1$. These are the following two cases. Empty sums and empty collections of additional bridge joins are omitted.

{\noindent\bf Case 1:} $r_{t-1}=r_t=2$ and $r_i=1$ for $i\in[t-2]$.

Let $a,b\in V(C_{t-1})$ and $c,d\in V(C_t)$ be the neighbors of $u$ in these components, and let
\[
 H=G[V(C_{t-1})\cup V(C_t)\cup\{u\}].
\]
For every $i\in[t-2]$, let $w_i$ be the unique neighbor of $u$ in $C_i$, and choose a maximum dissociation set $D_i$ of $C_i$ avoiding $w_i$. Any dissociation set of $H$ can be combined with these $D_i$ to form a dissociation set of $G$. By Claim 1,
\begin{align*}
 {\rm diss}(H)+\sum_{i=1}^{t-2}{\rm diss}(C_i)
 &\le {\rm diss}(G)\\
 &={\rm diss}(G-u)=\sum_{i=1}^{t}{\rm diss}(C_i).
\end{align*}
Since $C_{t-1}\cup C_t=H-u$ is an induced subgraph of $H$, we obtain
\[
 {\rm diss}(H)={\rm diss}(C_{t-1})+{\rm diss}(C_t).
\]
Consequently, for every maximum dissociation set $D'$ of $C_{t-1}\cup C_t$, either $|D'\cap\{a,b,c,d\}|\ge2$, or $D'\cap\{a,b,c,d\}=\{z\}$ and $d_{(C_{t-1}\cup C_t)[D']}(z)=1$. Otherwise, $D'\cup\{u\}$ would be a dissociation set of $H$ with more than ${\rm diss}(H)$ vertices. This proves condition (1) of $\bigoplus_6$.

Claim 1 also shows that no induced cycle of length $1$ modulo $3$ contains $u$. An induced $a$--$b$ path in $C_{t-1}$ or an induced $c$--$d$ path in $C_t$ of length $2$ modulo $3$, together with $u$, would yield such a cycle. Hence condition (2) of $\bigoplus_6$ holds. We may therefore construct $H=C_{t-1}\bigoplus_6 C_t$ and then recover
\[
 G=H\biguplus C_1\biguplus\cdots\biguplus C_{t-2}
\]
by adding the bridges $uw_i$. Thus $G\in\mathcal{G}$.

{\noindent\bf Case 2:} $r_t=3$ and $r_i=1$ for $i\in[t-1]$.

Let $x,y,z$ be the three neighbors of $u$ in $C_t$, put $S=\{x,y,z\}$, and let $H=G[V(C_t)\cup\{u\}]$. For every $i\in[t-1]$, choose a maximum dissociation set $D_i$ of $C_i$ avoiding the unique neighbor $w_i$ of $u$ in $C_i$. As in Case 1, any dissociation set of $H$ can be combined with these sets, so
\[
 {\rm diss}(H)+\sum_{i=1}^{t-1}{\rm diss}(C_i)
 \le {\rm diss}(G)=\sum_{i=1}^{t}{\rm diss}(C_i).
\]
The reverse inequality ${\rm diss}(H)\ge {\rm diss}(C_t)$ is immediate, and hence ${\rm diss}(H)={\rm diss}(C_t)$. It follows that, for every maximum dissociation set $D'$ of $C_t$, either $|D'\cap S|\ge2$, or $D'\cap S=\{v\}$ and $d_{C_t[D']}(v)=1$; otherwise $D'\cup\{u\}$ would be a larger dissociation set of $H$. Thus condition (1) of $\bigoplus_5$ holds.

If an induced path $P$ of length $2$ modulo $3$ had endpoints $s_1,s_2\in S$ and satisfied $V(P)\cap S=\{s_1,s_2\}$, then $P$ together with $u$ would be an induced cycle of length $1$ modulo $3$ in $G$. This contradicts Claim 1. Therefore condition (2) of $\bigoplus_5$ also holds, and $H=C_t\bigoplus_5\{u\}$. Joining the remaining components by their bridges gives
\begin{center}
\(
 G=(C_t\bigoplus_5\{u\})\biguplus C_1\biguplus\cdots\biguplus C_{t-1}
 \in\mathcal{G}.
\)
\end{center}
This completes the induction. $\qed$

\section*{Conclusion}

In this paper, we characterize all extremal graphs achieving equality in Theorem \ref{theorem1}. Since the dissociation number ${\rm diss}(G)$ is integer-valued, the lower bound $n-\frac{1}{3}(m+k+c_1)$ can be strengthened to $\lceil n-\frac{1}{3}(m+k+c_1)\rceil$. Therefore, it is also interesting to characterize all extremal graphs that satisfy ${\rm diss}(G)=\lceil n-\frac{1}{3}(m+k+c_1)\rceil$. The corresponding equality case for trees has been characterized by Li et al. \cite{li}.



\end{spacing}

\end{document}